\documentclass[11pt,amsart,leqno]{amsart}
\usepackage{amsmath}
\usepackage{amsfonts}
\usepackage{amssymb}
\usepackage{amsthm}
\usepackage{graphicx}
\usepackage[margin=2.5cm]{geometry}

\newtheorem{Proposition}{Proposition}

\theoremstyle{remark}
\title{A note on a result of Saks and Zygmund on additive functions of rectangles}

\author{Julià Cufí and Juan J. Donaire}

\thanks{Second author is supported in part by the Generalitat de Catalunya (grant 2017 SGR 395), the Spanish Ministerio de Ciencia e Innovaci\'on (projects  MTM2017-85666-P and PID2021-123151NB-I00) and the María de Maeztu Award.}
\address{Departament de Matem\`atiques, Universitat Aut\`onoma de Barcelona, 08193 Bellaterra, Barcelona, Catalonia\newline
Departament de Matem\`atiques, Universitat Aut\`onoma de Barcelona  and Centre de Recerca Matemàtica 08193 Bellaterra, Barcelona, Catalonia}
\email{julia.cufi@uab.cat}
\email{juanjesus.donaire@uab.cat}
\date{}

\begin{document}

\maketitle

\begin{abstract}
\noindent We modify the proof of the basic lemma of a paper of Saks and Zygmund on additive functions of rectangles and we weaken their hypotheses.
\end{abstract}

\section{Introduction.}\label{sec1}
 
 In the paper \cite{SZ}, Saks and Zygmund present a theorem on additive functions on rectangles that implies a slightly more general form of two theorems of Besicovitch about the set of removable singularities of a continuous or bounded analytic function on a simply connected domain.
 
 The main result (Theorem 5.2) is based on a lemma (Lemma 5.1) whose proof needs to be modified. At the same time their hypotheses can be slightly weakened.
 
\section{Discussion on the proof of Saks-Zygmund's lemma.}\label{sec1}

First of all let us introduce some definitions following \cite{SZ}.

We shall only consider rectangles and squares in the plane with sides parallel to the axis. A function $F$ of rectangles is said to be additive if $F(I_1\cup I_2)=F(I_1)+F(I_2)$, for any pair of adjacent rectangles $I_1$ and $I_2$.

For a rectangle $I$, we will denote by $\delta(I)$ its diameter and by $|I|$ its measure. The function of rectangles $F$ is said to be continuous if $F(I)\to 0$ as $\delta(I)\to 0$.

We consider the set of dyadic squares of the plane (meshes in \cite{SZ}). The dyadic squares of order $n$ of size $2^{-n}$ are the set of squares into which the plane is divided by the two systems of parallel lines
$$x=k\cdot 2^{-n}\;\;\;\; y=m\cdot 2^{-n},\;\;\; k,m\in\mathbb{Z}.
$$
Let $F$ be an additive and continuous function of rectangles, let $I_0$ be a fixed rectangle and $0\leq\alpha\leq 2$.

Let us define for $x\in I_0$,
$$\underline{F}_\alpha(x)=\liminf_{\delta(Q)\to 0}\frac{F(Q)}{\delta(Q)^\alpha}\;\;\text{ and }\;\;\underline{F}(x)=\liminf_{\delta(Q)\to 0}\frac{F(Q)}{|Q|},$$
where the limit is taken over arbitrary squares containing the point $x$.

For a set $E$ we will denote by $\Lambda_\alpha(E)$ the length of order $\alpha$ of $E$, that is to say its $H^\alpha$-Hausdorff measure. A set $E$ will be named a $B_\alpha$-set if $E$ is the union of a sequence of sets of finite length of order $\alpha$.

 Lemma 5.1 in \cite{SZ} states that if $F$ satisfies the conditions $\underline{F}_\alpha(x)\geq 0$ for $x\in I_0$ and $\underline{F}(x)\geq 0$ everywhere  in $I_0$ except for $x$ belonging to a $B_\alpha$-set, then $F(I_0)\geq 0$.

Taking into account the continuity of $F$, the argument is based on the reduction to the case in which $I_0$ is a dyadic square. In other words, one assumes that if an additive and continuous function of rectangles $F$ satisfies $F(Q_d)\geq 0$ for any dyadic square $Q_d$, then $F(I_0)\geq 0$ for every rectangle $I_0$. But this is false as the following example shows.

\begin{Proposition}\label{prop}
There is an additive and continuous function of rectangles $F$ such that $F(Q_d)\geq 0$ for every dyadic square $Q_d$ but $F(I_0)<0$ for some rectangle $I_0$.
\end{Proposition} 

\emph{Proof.}

A standard way to define an additive function $F$ of rectangles is the following one. Take a function $f$ defined on $\mathbb{R}^2$ and then, for any rectangle 
$$I=[x_1,x_2]\times [y_1,y_2]\;\; x_j,y_j\in\mathbb{R},$$
set $F(I)=f(x_2,y_2)+f(x_1,y_1)-f(x_1,y_2)-f(x_2,y_1)$, which is clearly additive.

Consider now the function $f$ defined as
$$f(x,y)=\left\{\begin{array}{ccl}
    1 &\text{ if }& y\not\in\mathbb{Q},  \\
    xy & \text{ if } &y\in\mathbb{Q}.
\end{array}\right.$$
Then, for any dyadic square $Q_d=[x_1,x_2]\times [y_1,y_2]$ one gets easily $F(Q_d)=(x_2-x_1)(y_2-y_1)=|Q_d|>0$, but taking $I=[x_1,x_2]\times [y_1,y_2]$ with $y_1\in\mathbb{Q}$, $y_1>0$ and $y_2\not\in\mathbb{Q}$, one gets $F(I_0)=(x_1-x_2)y_1<0$.

Note that this function $F$ is continuous since $\delta(I)\to 0$ implies that $x_2-x_1\to 0$ and $y_2-y_1\to 0$ if $I=[x_1,x_2]\times [y_1,y_2]$ and consequently, $F(I)\to 0$ in any case. \qed

\medskip

\emph{Remark.}

One could define another type of continuity for an additive function of rectangles $F$, namely that $F$ is continuous if $F(I)\to 0$ when $|I|\to 0$. This is a stronger notion of continuity than the one used in the quoted paper.

If one assumes this kind of continuity for $F$ then it is true that $F(Q_d)\geq 0$ for any dyadic square implies that $F(I_0)\geq 0$ for any rectangle $I_0$. This can be checked just by approximating $I_0$ by dyadic squares inside.

Note that the function $F$ exhibited in the proof of Proposition \ref{prop} is not continuous in this stronger sense.

\section{Modification of the proof of the lemma.}

In order to provide a proof of Lemma 5.1 in \cite{SZ} we will first show that the hypotheses in this Lemma imply that $F(Q)\geq 0$ for any square $Q\subset I_0$. Moreover we will remark that the hypotheses can be slightly weakened by assuming only that $F(I)$ satisfies the condition $\underline{F}_\alpha(x)\geq 0$ on $I_0\setminus\Sigma$, where $\Sigma$ is a set with $\Lambda_\alpha(\Sigma)=0$.

So we can state the following

\textbf{Lemma 5.1'.\;}\sl If an additive and continuous function of rectangles $F$ satisfies $\underline{F}_\alpha(x)\geq 0$, where $0\leq\alpha\leq 2$, for $x\in I_0\setminus\Sigma$, with $I_0$ a rectangle and $\Sigma\subset I_0$ with $\Lambda_\alpha(\Sigma)=0$ and $\underline{F}(x)\geq 0$ everywhere in $I_0$ except for x belonging to a $B_\alpha$-set, then $F(Q)\geq 0$ for every square $Q\subset I_0$.\rm

For the proof, fix a square $Q\subset I_0$ of sidelength $\ell$ and consider as dyadic squares the squares obtained by doing a dyadic division of $Q$. Consequently, at the n-th step one obtains $2^{2n}$ squares of sidelength $\ell/2^n$. Now Lemma 4.1 in \cite{SZ} for this familiy of dyadic squares holds (going to the unit square by a dilation and a translation) with assertion (i) with a different constant from the constant 32 appearing there.

Then one can follow the proof of Lemma 5.1 in \cite{SZ} where now $D$ is replaced by $D\cap Q$ and $D_i$ by $D_i\cap Q$. To see that the hypotheses of Lemma 5.1 can be weakened just note that we will have, following the notation in the proof of this Lemma,
$$I_0=\bigcup_{n=1}^\infty R_{i,n}\cup\Sigma,\text{ for }i=1,2,\cdots$$
and $\Lambda_\alpha(D_i)=\sum_{n=1}^\infty\Lambda_\alpha(D_{i,n})+\Lambda_\alpha(D_i\cap\Sigma)=\sum_{n=1}^\infty\Lambda_\alpha(D_{i,n})$ as in (5.2) in \cite{SZ}.

At the end for $G(I)=F(I)+\varepsilon |I|$ we will obtain 
$$F(Q)=G(Q)-\varepsilon |Q\geq -(C+|Q|)\varepsilon$$
for some constant $C$ and $\varepsilon$ arbitrarily small, that gives $F(Q)\geq 0$.

Therefore Lemma 5.1 in \cite{SZ} will be obtained with the weak hypothesis $\underline{F}_\alpha(x)\geq 0$ for $x\in I_0\setminus\Sigma$ with $\Lambda_\alpha(\Sigma)=0$ if we stablish the following fact.

\begin{Proposition}
Let $F$ be an additive and continuous function of rectangles such that $F(Q)\geq 0$ for any square $Q\subset I_0$, $I_0$ some rectangle. Then one has $F(I_0)\geq 0$.
\end{Proposition}

\emph{Proof.}

Assume that $I_0$ has sides of lengths $\ell_0>\ell_1$. Put inside $I_0$ adjacent squares of sidelength $\ell_1$ as many times as possible and let $R_1$ the remaining rectangle.

Following this process with the rectangle $R_1$, after a number of iterations we will obtain a decomposition
$$I_0=Q_1\cup Q_2\cup\cdots \cup Q_n\cup R_n$$
with $Q_1, Q_2,\cdots ,Q_n$ squares and $R_n$ a rectangle with largest side $\ell_n$. If we show that $\ell_n\simeq\delta(R_n)\to 0$ as $n\to\infty$ we are done.
  \begin{center}
	\includegraphics[height=4.5cm]{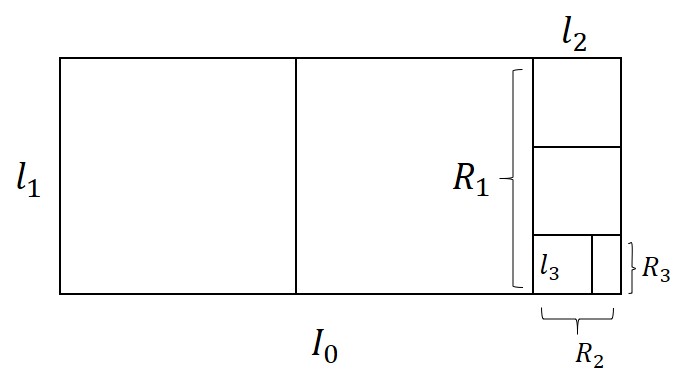}
\end{center}

If $\ell_{n+1}\leq \ell_n/2$ for any $n$, the result is clear. However it can occur that for some $n$ one has $\ell_{n+1}> \ell_n/2$, but in this case, $\ell_{n+2}\leq \ell_n/2$. In fact, $\ell_{n+1}$ only fits once in $\ell_n$ and one has 
$$\ell_{n+2}=\ell_{n}-\ell_{n+1}\leq \ell_n-\frac{\ell_n}{2}=\frac{\ell_n}{2}.$$
So, at every step or at every two steps, the sidelength is reduced at least by a factor $1/2$ and, consequently $\ell_n\to 0$, since $(\ell_n)$ is decreasing.\qed

\section{Modified results}

The comments we have made above imply that the main results in \cite{SZ} can be stated with less restrictive hypotheses.

First of all in Theorem 5.1 we can assume that the function $F(I)$ satisfies the condition $\underline{F}_\alpha (x)\geq 0$ in a rectangle $I_0$ except for a set $\Sigma$ with $\Lambda_\alpha(\Sigma)=0$, because then the difference function $\Delta(I)=F(I)-\Psi(I)$ will satisfy $\underline{\Delta}_\alpha(x)\geq 0$ on $I_0\setminus\Sigma$ and Lemma 5.1' will give $\Delta(I)\geq 0$, for $I\subset I_0$.

As a consequence, Theorem 5.2 will be true assuming that the property $(\ell_\alpha)$ that means
$$F'_\alpha(x):=\lim_{\delta(Q)\to 0}\frac{F(Q)}{\delta(Q)}=0,$$
where $Q$ runs on squares containing the point $x$, is true for $x\in I_0\setminus\Sigma$ with $\Lambda_\alpha(\Sigma)=0$.

Finally, the result in section (6) of \cite{SZ} can be stated with the weaker hypothesis that the complex function $f(z)$ is continuous out of a set $\Sigma$ with $\Lambda_1(\Sigma)=0$. The reason is that then the interval functions $U(I)$ and $V(I)$ defined from
$$\int_{\partial I} f(z)\, dz=U(I)+iV(I)$$
will satisfy the condition $(\ell_1)$, that is to say
$$\lim\limits_{{\delta(Q)\to 0}\atop{x\in Q}}\frac{U(Q)}{\delta (Q)}=\lim\limits_{{\delta(Q)\to 0}\atop{x\in Q}}\frac{V(Q)}{\delta (Q)}=0$$
at the points $x$ not in $\Sigma$ and we can apply Theorem 5.2 with the modified hypotheses.

\end{document}